\documentclass[a4paper]{article}
\usepackage{amsmath,amsfonts,amssymb,amsthm}
\usepackage{latexsym, xspace, enumerate}
\usepackage[mathscr]{eucal}
\usepackage[all]{xy}
\usepackage{hyperref}
\usepackage{amscd}
\usepackage{mathtools}
\usepackage{vmargin}
\usepackage{rotating}

\setmarginsrb{20mm}{10mm}{20mm}{25mm}
         {10mm}{7mm}{10mm}{10mm}

\newtheorem{theorem}{Theorem}[section]
\newtheorem*{corollary1}{Lemma 1.3}
\newtheorem*{theorem2}{Theorem 1.1}
\newtheorem*{corollary3}{Corollary 1.2}
\newtheorem*{theorem4}{Theorem 1.5}
\newtheorem*{theorem5}{Theorem 1.7}
\newtheorem*{theorem6}{Theorem 1.9}
\newtheorem*{lemma6}{Proposition 1.11}

\newtheorem{proposition}[theorem]{Proposition}

\newtheorem{lemma}[theorem]{Lemma}

\newtheorem{corollary}[theorem]{Corollary}

\theoremstyle{definition}

\newtheorem{example}[theorem]{Example}
\newtheorem{remark}[theorem]{Remark}

\newcommand{\N}{\mathbb N}

\newcommand{\Z}{\mathbb Z}

\def\Sym{\mathrm{Sym}}
\newcommand{\id}{\text{id}}

\makeatletter
\def\Ddots{\mathinner{\mkern1mu\raise\p@
\vbox{\kern7\p@\hbox{.}}\mkern2mu
\raise4\p@\hbox{.}\mkern2mu\raise7\p@\hbox{.}\mkern1mu}}
\makeatother

\title{Residual properties of free products}
\author{Federico Berlai\footnote{Universit\"{a}t Wien, Fakult\"{a}t f\"{u}r Mathematik, Oskar-Morgenstern-Platz 1, 1090 Wien, Austria.
\newline
E-mail address: federico.berlai@univie.ac.at.}}
\date{}

\begin{document}
\maketitle
\begin{abstract}
Let $\mathcal{C}$ be a class of groups. We give sufficient conditions ensuring that a free product of residually $\mathcal{C}$ groups is again residually $\mathcal{C}$, and
analogous conditions are given for LE-$\mathcal{C}$ groups. As a corollary, we obtain that the class of residually amenable groups and the one of 
LEA groups are closed under taking free products.

Moreover, we consider the pro-$\mathcal{C}$ topology and we characterize special HNN extensions and amalgamated free products that are residually $\mathcal{C}$,
where $\mathcal{C}$ is a suitable class of groups. In this way, we describe special HNN extensions and amalgamated free products that are residually amenable.
\end{abstract}

\bigskip
\noindent 
\textrm{\small{Key words: root property, residually amenable group, sofic group, free product, pro-$\mathcal{C}$ topology, proamenable topology.}}

\section{Introduction}

Amenable groups were defined in 1929 by von Neumann in his work on the Banach-Tarski paradox and, since then, the study of amenability has been remarkably fruitful
in a variety of branches of mathematics. Recently, weak variants of amenability have attracted considerable interest. Examples are the Haagerup property,
also known as Gromov's a-T-menability, and coarse amenability, also termed as Guoliang Yu's Property A. 
These notions appeared in the study of analytic properties of infinite groups.

The main focus of this paper are other, more algebraic, weak forms of amenability: the notions of residual amenability and of LEA (Locally 
Embeddable into Amenable) groups.
Residual amenability is a common generalization of amenability and of residual finiteness.
While the stability, or permanence, properties of both the classes of amenable groups and of residually finite groups
are well known, the stability properties of the class of residually amenable groups have not yet been investigated. 
See, for instance, \cite{Sch01} for corresponding open questions.

In \cite{CSC} residually amenable and LEA groups
appear as fundamental examples of sofic \cite{Gro,Weiss} and hyperlinear \cite{Pes08} groups. 
Thus, on the one hand, it seems reasonable to study the (non-)permanence properties of both the classes of residually amenable and LEA groups, 
in order to have a better understanding also of the latter classes of sofic and hyperlinear groups.

On the other hand, after more than ten years since sofic groups were defined, the following question is still open: 
Is there a non-sofic/non-hyperlinear group?
There are few examples of groups which are known to be sofic but not residually amenable. Hence, new examples of non-residually amenable
groups are of great value for the above question, and one way to produce these new examples is to study the non-permanence properties of the
class of residually amenable groups.

On yet another perspective, residually amenable groups have recently been investigated in the context of the theory of $L^2$-invariants. In \cite{Clair97} the approximation of $L^2$-Betti numbers
is studied, and in \cite{Lueck} the connections between residually amenable groups and the $L^2$-invariants
theory are explained, along with deep conjectures they satisfy, such as the Approximation Conjecture on the $L^2$-Betti numbers, 
and the Determinant Conjecture \cite[Conjectures 13.1 and 13.2]{Lueck}. 
Therefore, the (non-)stability properties of the class of residually amenable groups will also yield new examples in this context.

\medskip

The first aim of the paper is to show that the class of residually amenable groups is closed under free products. To this end, 
in Section \ref{sec2} we consider the classes of residually amenable and residually elementary amenable groups. We show that the classical
result of Gruenberg \cite[Theorem 4.1]{Gruenberg57} cannot be applied to these classes, that is, amenability and elementary amenability are not 
root classes (see next section for definitions). 

In particular, a new method is required to deal with free products of residually amenable and residually elementary amenable groups.
This is the focus of Section \ref{sec3}, where we prove that the classes of residually $\mathcal{C}$ groups are closed under free products, when
$\mathcal{C}$ satisfies certain conditions. From this, we deduce that the class of residually amenable groups is indeed closed under free products.
After that, we obtain analogous results for free products of LEA groups and LE-$\mathcal{C}$ groups.

In Section \ref{section-proamenable}, we consider the pro-$\mathcal{C}$ topology and the pro-$\mathcal{C}$ completion of a group, and
we study special HNN extensions and amalgamated free products which are residually $\mathcal{C}$. As a corollary, we characterize particular HNN extensions and amalgamated free products which are 
residually amenable. Our result is motivated by recent developments
in the theory of right-angled Artin groups, and by the analogous fact known for residually finite groups.

\subsection{Background}

We start with standard terminology and known facts that are used in this paper.

A \emph{class} of groups is a family of groups closed under isomorphic images. Let $\mathcal{C}$ be a class of groups. Then a group $G$ is
\emph{residually} $\mathcal{C}$ if for every non-trivial element $g\in G$ there exist a group $C\in\mathcal{C}$ and a surjective homomorphism
$\varphi\colon G\twoheadrightarrow C$ such that $\varphi(g)\neq e_C$. Moreover, $G$ is \emph{fully residually} $\mathcal{C}$ if for each finite subset $H\subseteq G\setminus\{e_G\}$ there
exist a group $C\in\mathcal{C}$ and a surjective homomorphism $\varphi\colon G\twoheadrightarrow C$ such that $e_C\notin \varphi(H)$. 
If $\mathcal{C}$ is closed under finite direct products, then a group is residually $\mathcal{C}$ if and only if it is fully 
residually~$\mathcal{C}$. 

A related notion is the one of local embeddability.
Let $G$ and $C$ be two groups and $K\subseteq G$ a finite subset. Then
a map $\varphi\colon G\to C$ is called a $K$-\emph{almost}-\emph{homomorphism} if it satisfies the following conditions:
\begin{itemize}
 \item[$(a1)$] $\varphi(k_1k_2)=\varphi(k_1)\varphi(k_2)$ for all $k_1$, $k_2\in K$;
 \item[$(a2)$] $\varphi\restriction_K$ is injective.
\end{itemize}
Let $\mathcal{C}$ be a class of groups. Then a group $G$ is \emph{locally embeddable into the class} $\mathcal{C}$ ($G$ is a LE-$\mathcal{C}$ group) if for each finite 
subset $K\subseteq G$ there exist a group $C\in\mathcal{C}$ and a $K$-almost-homomorphism $\varphi\colon G\to C$. 
For instance, if $\mathcal{C}$ is the class of amenable groups then we get the notion of LEA groups, also known as Gromov's initially subamenable 
groups \cite[page 133]{Gro}.
If $\mathcal{C}$ is the class of finite groups, we obtain LEF groups, see e.g. \cite{VeGo97}.

If the class $\mathcal{C}$ is closed under finite direct products, then being residually $\mathcal{C}$ implies being LE-$\mathcal{C}$ 
\cite[Corollary~7.1.14]{CSC}. Moreover, if
the class $\mathcal{C}$ is also closed under taking subgroups, then a finitely presented group is residually $\mathcal{C}$ if and only if it is LE-$\mathcal{C}$ \cite[Corollary 7.1.22]{CSC}.

The study of residual properties began with residual finiteness: in the 1930s Levi proved that free groups are residually finite, and soon after
that more proofs appeared, as the one of Hall~\cite{Hall49}.
Remarkable was the work that Mal'cev started in the 1940s on linear groups and residual finiteness. He proved that finitely generated linear groups are residually finite, and 
that finitely generated residually finite groups are Hopfian \cite{Mal}. He also showed that a split extension of a finitely generated residually finite group 
by a residually finite
group is residually finite. Recently, this result of Mal'cev concerning split extensions was generalized by Arzhantseva and Gal \cite[Theorem 7]{ArGal}.

Since the 1960s, the interest in residual properties evolved into the study of residually nilpotent, residually free, fully residually free 
(also known as limit groups or $\omega$-residually free groups) and residually solvable groups. 
Residually amenable groups are a common generalization of all these classes.

In 1957 Gruenberg introduced the concept of \emph{root property}, or \emph{root class} of groups, 
namely, a non-trivial class $\mathcal{C}$ such that
\begin{enumerate}
 \item $\mathcal{C}$ is closed under taking subgroups;
 \item $\mathcal{C}$ is closed under taking finite direct products;
 \item \emph{Gruenberg Condition}: for any chain $K\unlhd H\unlhd G$ such that $G/H$ and $H/K$ are in $\mathcal{C}$, 
 there exists $L\unlhd G$ such that $L$ is contained in $K$ and $G/L\in\mathcal{C}$.
\end{enumerate}
Finite groups, finite $p$-groups and solvable groups are known to satisfy the Gruenberg condition \cite{Gruenberg57}, while nilpotent groups do not 
satisfy it.
Using the fact that finite groups form a root class,
Gruenberg proved that the free product of residually finite groups is residually finite \cite[Theorem 4.1]{Gruenberg57}.

\medskip

Before stating our results, let us fix some notations. If $\mathcal{A}$ and $\mathcal{B}$ are two properties of groups, 
then a group $G$ is $\mathcal{A}$-by-$\mathcal{B}$ if it admits a normal subgroup $N$ satisfying $\mathcal{A}$ such that the quotient 
group $G/N$ satisfies $\mathcal{B}$.

The identity of a group $G$ is denoted by $e_G$, or simply by $e$ if the group is clear from the context.
The conjugate $gHg^{-1}$ of a subgroup $H\leq G$ is
denoted by $H^g$, and the free group on $n$ free generators is denoted by $F_n$.
If $X$ is a set, $\Sym(X)$ denotes the group of permutations of $X$, $\Sym_0(X)$ the group of permutations of $X$ with finite support and 
$\Sym_0^+(X)$ the alternating subgroup of $\Sym_0(X)$.

\subsection{Main results}

Our first result is the following generalization of the classical theorem of Gruenberg, 
who showed that free products of residually finite groups are residually finite \cite[Theorem~4.1]{Gruenberg57}.
\begin{theorem}\label{theoremResi}
Let $\mathcal{C}$ be a class of groups that contains a root class $\mathcal{R}$. Suppose that
\begin{enumerate}
 \item $\mathcal{C}$ is closed under finite direct products;
 \item every $\mathcal{R}$-by-$\,\mathcal{C}$ group belongs to $\mathcal{C}$;
 \item for every group in $\mathcal{C}$ there is a group in $\mathcal{R}$ of the same cardinality.
\end{enumerate} 
Then the class of residually $\mathcal{C}$ groups is closed under free products.
\end{theorem}
We deduce that the class of residually amenable groups is closed under free products. This answers an open question
raised by Schick \cite{Sch01}.
\begin{corollary}\label{corollary.resam}
The classes of residually amenable groups and of residually elementary amenable groups are closed under free products.
\end{corollary}

To prove Corollary \ref{corollary.resam} from Theorem \ref{theoremResi}, we use the following result.
\begin{lemma}\label{RS-by-A}
If a group is (residually solvable)-by-amenable then it is residually amenable. 

If a group is (residually solvable)-by-(elementary amenable) then it is residually elementary amenable. 
\end{lemma}

The next fact is of independent interest. See Example \ref{ex-amenableNOGRUENBERG} for the proof.
\begin{proposition}\label{RA-by-A}
The class of residually amenable (resp. residually elementary amenable) groups is not closed under extensions with amenable (resp. elementary 
amenable) groups.
\end{proposition}

Our second main result is the following theorem on free products of LE-$\mathcal{C}$ groups.
\begin{theorem}\label{theoremLocalEmbeddings}
Let $\mathcal{C}$ be a class of groups, and suppose that
\begin{itemize}
 \item[$(a)$] $\mathcal{C}$ is closed under free products
\end{itemize} 

or
\begin{itemize}
 \item[$(b)$] $\mathcal{C}$ is closed under finite direct products and the free product of residually $\mathcal{C}$ groups is residually 
 $\mathcal{C}$.
\end{itemize}
Then the class of LE-$\mathcal{C}$ groups is closed under free products.
\end{theorem}
\begin{corollary}\label{LEA-free}
Let $\mathcal{C}$ be one of the following classes:
\begin{enumerate}
 \item finite groups;
 \item finite $p$-groups;
 \item solvable groups;
 \item free groups;
 \item elementary amenable groups;
 \item amenable groups.
\end{enumerate}
Then the class of LE-$\mathcal{C}$ groups is closed under free products.
\end{corollary}

After recalling the definition of pro-$\mathcal{C}$ topology, we focus our attention on special HNN extensions and doubles of residually amenable groups 
(see Section \ref{section-proamenable} for the definitions).

\begin{theorem}\label{specialHNN}
Let $A$ be a group, $H\leq A$ and consider the special HNN extension $G:=A\ast_{\id}$ amalgamating the subgroup $H$. 
Let $\mathcal{C}$ be a class of groups that contains a root class $\mathcal{R}$. Suppose that
\begin{enumerate}
 \item $\mathcal{C}$ is closed under subgroups and finite direct products;
 \item every $\mathcal{R}$-by-$\,\mathcal{C}$ group belongs to $\mathcal{C}$;
 \item for every group in $\mathcal{C}$ there is a group in $\mathcal{R}$ of the same cardinality.
\end{enumerate}
Then $G$ is residually $\mathcal{C}$ if and only if $A$ is residually $\mathcal{C}$ and $H$ is closed in the pro-$\mathcal{C}$ topology of $A$.
\end{theorem}

\begin{corollary}\label{cor-HNNspecial}
A special HNN extension $G:=A\ast_{\id}$ is residually amenable if and only if $A$ is residually amenable and the amalgamated subgroup is closed 
in the proamenable topology of $A$.
\end{corollary}

\begin{theorem}\label{theo.bingobongo}
Let $A$ be a group, $H\leq A$ a subgroup and consider the double $G=A\ast_H A$. 
Let $\mathcal{C}$ be a class of groups that contains a root class $\mathcal{R}$. Suppose that
\begin{enumerate}
 \item $\mathcal{C}$ is closed under subgroups and finite direct products;
 \item every $\mathcal{R}$-by-$\,\mathcal{C}$ group belongs to $\mathcal{C}$;
 \item for every group in $\mathcal{C}$ there is a group in $\mathcal{R}$ of the same cardinality.
\end{enumerate}
Then $G$ is residually $\mathcal{C}$ if and only if $A$ is residually $\mathcal{C}$ and $H$ is closed in the pro-$\mathcal{C}$ topology of $A$.
\end{theorem}

\begin{corollary}\label{cor.double}
A double $A\ast_HA$ is residually amenable if and only if $A$ is residually amenable and $H$ is closed in the proamenable topology of $A$. 
\end{corollary}

Finally, we extend results on HNN extensions \cite[Lemma~1.2]{RaVa} to residual amenability.
\begin{proposition}\label{stabilityHNN}
Let $A$ be a group, $H$, $K$ two subgroups of $A$ and $\varphi\colon H\to K$ an isomorphism. Suppose that there exists an automorphism 
$\alpha\colon A\to A$ such that $\alpha\restriction_H=\varphi$.
If $A$ is amenable then the HNN extension $A\ast_{\varphi}$ is residually amenable. 
\end{proposition}

\section{The class of residually amenable groups}\label{sec2}

In this section, we recall some properties of the class of amenable and elementary amenable groups, and then we focus on residual amenability. 
Very basic permanence properties of this class of groups can be obtained using \cite[\S 2.2]{CSC}: the known stability properties
of residually finite groups proved in that book can be instantly extended to the corresponding stability properties of residual amenability. 
The only result of \cite[\S 2.2]{CSC} that is not
immediate to generalize is \cite[Lemma 2.2.11]{CSC}, on the residual finiteness of virtually residually finite groups.
We prove the analogous fact for residually amenable groups in Lemma \ref{virtually.residually.amenable}.

\smallskip\smallskip

Finite groups and solvable groups are amenable, and the class of discrete amenable groups is closed under taking subgroups, quotients, extensions, and direct limits.
If a group contains a copy of a free non-abelian group then it is not amenable, and a virtually amenable groups is amenable. We refer to \cite[\S 2.2]{CSC} for the proofs of all of 
these facts.

The class of elementary amenable groups is the smallest class that contains the finite and the solvable groups, and which is closed under the four operations 
just mentioned, namely, closed under taking subgroups, quotients, extensions, and direct limits.
Every elementary amenable group is amenable, but it is now well known that there are many amenable groups which are not elementary amenable, e.g. \cite{Gri84}.

The following fact is probably known to specialists, although we have not located an exact reference in the literature.

\begin{proposition}\label{prod-not-amenable}
The group $G:=\prod_{ \N}\Sym_0(\Z)$ contains copies of free non-abelian groups. In particular, $G$ is not amenable.

\begin{proof}

The group $\Sym_0(\Z)$ is elementary amenable because it is locally finite.
Let $p_i$ be the $(2i-1)$-th prime number and $q_i$ be the $2i$-th prime number, for all $i\geq 1$. 
The families $\{p_i\}_{i\geq1}$ and $\{q_i\}_{i\geq1}$ are disjoint and their union is exactly the
set of all prime numbers.
Define, for each $i\geq 1$, the permutations $g_i$ and $\gamma_i$ of $\Z$ with finite support
$$g_i:=(0,p_i,2p_i,\dots,ip_i),\qquad \gamma_i:=(0,q_i,2q_i,\dots,iq_i)$$
and consider the elements
\begin{equation*}\label{diabolic.elements}
g:=(g_n)_{n\in\N},\qquad \gamma:=(\gamma_n)_{n\in\N},
\end{equation*}
letting $g_0=\id_{\Z}=\gamma_0$. In particular, the $n$-th components $g_n$ and $\gamma_n$ of $g$ and $\gamma$ are cycles of length $n+1$.
Using the Ping-Pong Lemma we prove that 
\begin{equation*}\label{free-subgroup}
\langle g,\gamma\rangle\cong\langle g\rangle\ast \langle\gamma\rangle=F_2.
\end{equation*}
Consider the following sets:
$$X_{\pm}:=\bigl\{w\in G\mid w=g^{\pm n}w', n\in\N\setminus\{0\},\text{ $w'$ is reduced and does not begin with }g^i, 
\quad i\in\Z\setminus\{0\}\bigr\},$$
$$Y_{\pm}:=\bigl\{w\in G\mid w=\gamma^{\pm n}w', n\in\N\setminus\{0\},\text{ $w'$ is reduced and does not begin with }\gamma^i, 
\quad i\in\Z\setminus\{0\}\bigr\}.$$
These four sets are disjoint and $e_G$ does not belong to any of them, because $g$ and $\gamma$ are torsion-free elements of $G$.
As 
$$g(G\setminus X_-)\subseteq X_+,\qquad g^{-1}(G\setminus X_+)\subseteq X_-$$
and
$$\gamma(G\setminus Y_-)\subseteq Y_+,\qquad \gamma^{-1}(G\setminus Y_+)\subseteq Y_-,$$
we can apply the Ping-Pong Lemma to conclude that $\langle g,\gamma\rangle\cong F_2$, so the discrete group $G$ is not amenable.

\end{proof}
\end{proposition}

Thus, there exist amenable groups $G$ such that the power $\prod_{\N}G$ is not amenable.
However this cannot occur whenever $G$ is finite as in that case $\prod_{\N}G$ is locally finite.

\begin{example}\label{ex-amenableNOGRUENBERG}
Let $A:=\Sym_0(\Z)$, consider the unrestricted wreath product $G:=A\wr \Z=A^{\Z}\rtimes_\sigma\Z$ and let
\begin{equation}\label{enne}
 N=:\{(a_i)_{i\in\Z}\in A^\Z\mid a_0=e_A\}\unlhd A^\Z.
\end{equation}
Then
\begin{enumerate}
 \item the subnormal chain of subgroups $N\trianglelefteq A^\Z\trianglelefteq G$ does not satisfy the Gruenberg condition;
 \item the group $G$ is sofic but not residually amenable.
\end{enumerate}

\medskip

The groups $G/A^\Z\cong \Z$ and $A^\Z/N\cong A$ are both elementary amenable.
Suppose now that the Gruenberg condition holds, so that there exists a normal subgroup $L\unlhd G$ such that $L\subseteq N$ and $G/L$ is amenable. 

As $L\unlhd G$, we have that
\begin{equation*}
L\subseteq \bigcap_{\substack{g=(e,z)\\z\in \Z}}N^g.
\end{equation*} 
Let $g=(e,z)$, then
\begin{equation*}
\begin{split}
N^g&=\{(e,z)(n,0)(e,z)^{-1}\mid n\in N\}=\bigl\{\bigl(\sigma_z(n),z\bigr)(e,-z)\mid n\in N\bigr\}\\
&=\bigl\{\bigl(\sigma_z(n),0\bigr)\mid n\in N\bigr\}=\{(a_i)_{i\in\Z}\in A^\Z\mid a_z=e_A\},
\end{split}
\end{equation*}
thus $L$ is contained in
$$\bigcap_{z\in\Z}\{(a_i)_{i\in\Z}\in A^\Z\mid a_z=e_A\}=\{e_{A^\Z}\},$$
that is, $L=\{e_G\}$.
This implies that $G/L=G$ is amenable, and this is a contradiction of Proposition \ref{prod-not-amenable}.

\smallskip

Note that $G$ is (residually amenable)-by-amenable, thus it is sofic \cite[Proposition 7.5.14]{CSC};
more precisely it is elementary sofic \cite{Cor11}.

\smallskip

Let us prove that $G$ is not residually amenable. 
Consider the infinite simple group $B:=\Sym_0^+(\Z)\trianglelefteq A$. 
The quotient group $A^\Z/B^\Z\cong \Z_2^\Z$ 
is amenable because it is abelian.
This implies that $B^\Z$ is not amenable, because amenability is closed under extensions and $A^\Z$ is not amenable by Proposition \ref{prod-not-amenable}.

The group $B\wr \Z$ is a subgroup of $G=A\wr \Z$, so if $B\wr \Z$ is not 
residually amenable then neither is $G$.
Let $N\trianglelefteq B\wr\Z$ be a normal subgroup with amenable quotient, then the normal subgroup $N\cap B^\Z\trianglelefteq B\wr \Z$ 
is not trivial. Let $e_G\neq(n,0)\in N\cap B^\Z$ with $n_i\neq e_B$ for some $i\in\Z$. 
The center of $B$ is trivial because $B$ is simple and non-abelian, so there exists $b\in B$ such that $a:=[b,n_i]\neq e_B$. 

Consider the new element $\beta\in B^\Z$
\begin{equation*}\beta_j=\begin{cases}e_B&\quad\text{if }j\neq i,\\ b\neq e_B&\quad\text{if }j=i,\end{cases}\end{equation*}
then $\alpha:=[n,\beta]\in B^\Z$ is an element with only one non-trivial component:
\begin{equation*}\alpha_j=\begin{cases}e_B&\quad\text{if }j\neq i,\\ a\neq e_B&\quad\text{if }j=i.\end{cases}\end{equation*}
An easy computation shows that
$$(g,z)(\alpha,0)(g,z)^{-1}=(g\sigma_z(\alpha)g^{-1},0)\qquad \forall (g,z)\in B\wr\Z,$$
and these elements are all elements of $N\cap B^\Z$, because $N\cap B^\Z\unlhd B\wr\Z$.

As $e_B\neq a\in B$ and $B$ is simple, we have that
$$\{hah^{-1}\mid h\in B\}=B$$
and that
\begin{equation*}\begin{split}
B^\Z&\supseteq N\cap B^\Z \supseteq\bigl\{(g,z)(\alpha,0)(g,z)^{-1}\mid (g,z)\in B\wr\Z \bigr\}\\&
=\bigl\{(g\sigma_z(\alpha)g^{-1},0)\mid (g,z)\in B\wr\Z\bigr\}=\bigl\{\sigma_z(g\alpha g^{-1})\mid (g,z)\in B\wr\Z\bigr\}=B^\Z,
\end{split}\end{equation*}
that is, $N\cap B^\Z=B^\Z$.

This means that $B\wr \Z$ is not residually amenable, because for the elements $(b,0)\in B\wr \Z$ there is no normal subgroup $N$ with amenable
quotient, such that $(b,0)\notin N$. Hence the group
$$A\wr \Z= \Bigl(\prod_{z\in\Z}\Sym_0(\Z)\Bigr)\rtimes\Z$$
is not residually amenable. 
\end{example}

Proposition \ref{RA-by-A} follows from what we just proved.

\begin{lemma}
\label{virtually.residually.amenable}
Let $G$ be a virtually residually amenable group. Then $G$ is residually amenable.
\begin{proof}
Let $H\leq G$ be a residually amenable subgroup of finite index of $G$. 
Then there exists $N\unlhd G$ of finite index in $G$ such that $N\leq H$. Hence
$N$ is residually amenable.

Let $g\in G\setminus\{e\}$ be a non-trivial element, we want to find a normal subgroup of $G$ not containing $g$ such that the quotient is amenable.
If $g\in G\setminus N$ then it is sufficient to consider the finite quotient~$G/N$. 

Suppose that $g\in N\setminus\{e\}$.
As $N$ is residually amenable, there exists $K\unlhd N$ such that $g\notin K$ and $N/K$ is amenable.
Let
$$G/N=\{g_1N,g_2N,\dots,g_rN\}$$
and define $L:=\bigcap_{i=1}^rg_iKg_i^{-1}$. Then $L\unlhd G$, $g\notin L$ and $N/L$ has finite index in $G/L$.

We now show that $N/L$ is amenable, that is, $G/L$ is virtually amenable.
In fact, as $N\unlhd G$, every conjugate $g_iKg_i^{-1}$ is again a normal subgroup of $N$, so
$$N/L\leq \bigl(N/g_1Kg_1^{-1}\bigr)\times  \dots\times \bigl(N/g_rKg_r^{-1}\bigr)\cong \bigl( N/K\bigr)^r,$$
Thus $N/L$ is amenable because it is a subgroup of an amenable group.

Hence $G/L$ is indeed virtually amenable, that is, amenable. This concludes the proof, as we have shown that for every non-trivial element $g$ 
there exists $L\trianglelefteq G$ which does not contain $g$ and with amenable quotient.
 
\end{proof}
\end{lemma}

\section{Proof of the main theorems}\label{sec3}

We saw in the previous section that amenability is not a root property, thus \cite[Theorem 4.1]{Gruenberg57} cannot apply to show that 
residual amenability is preserved under free products.

In the following lemma we recall some easy properties implied by the Gruenberg condition. 
The second claim is false in general if the class $\mathcal{C}$ does not satisfy the Gruenberg condition. Proposition \ref{RA-by-A} and Example 
\ref{ex-amenableNOGRUENBERG} provide a counterexample.

\begin{lemma}\label{folklore}
If a class $\mathcal{C}$ satisfies the Gruenberg condition then 
\begin{enumerate}
 \item $\mathcal{C}$ is closed under extensions;
 \item a group which is (residually $\mathcal{C}$)-by-$\,\mathcal{C}$ is residually $\mathcal{C}$ \cite[Lemma 1.5]{Gruenberg57};
 \item free groups are residually $\mathcal{C}$.
\end{enumerate}
\begin{proof}
\begin{enumerate}
 \item 
Consider a group $G$ with a normal subgroup $N$ such that $N$ and $G/N$ are in $\mathcal{C}$. Apply the Gruenberg condition to the subnormal
chain
$\{e\}\unlhd N\unlhd G$
to conclude that $G\in\mathcal{C}$.

\item

Let $G$ be a (residually $\mathcal{C}$)-by-$\,\mathcal{C}$ group, with a normal residually $\mathcal{C}$ subgroup $N$ such that $G/N\in 
\mathcal{C}$.
Consider an element $g\in G\setminus\{e_G\}$, if $g\notin N$ then the canonical projection $\pi\colon G\twoheadrightarrow G/N$ does not map $g$ to the
identity of $G/N$. If $g\in N$ then there exists a surjective homomorphism $\varphi\colon N\twoheadrightarrow C$ such that $C\in \mathcal{C}$ and $\varphi(g)\neq e_C$, 
that is, $g\notin K=\ker \varphi$.
As $\varphi$ is surjective, we have that $N/K\cong C\in \mathcal{C}$, so in the subnormal chain
$K\trianglelefteq N\trianglelefteq G$
both quotients are groups in $\mathcal{C}$. Thus, applying the Gruenberg condition, we obtain a subgroup $L\trianglelefteq G$ contained in $K$ such
that $G/L\in\mathcal{C}$. As $L\subseteq K$, we have that $g\notin L$, and the canonical projection $\psi\colon G\twoheadrightarrow G/L$ does not map $g$ to
the identity of $G/L$. This means that $G$ is residually~$\mathcal{C}$.

\item

Let $G\in\mathcal{C}$ be a non-trivial group, and $g$ be a non-trivial element of $G$. Then either $g$ has infinite order or finite order.
In the first case, $\mathcal{C}$ contains all finitely generated torsion-free nilpotent groups. In the second, $\mathcal{C}$ contains all finite $p$-groups,
for $p$ a prime dividing the order of $g$.

Free groups are both residually finitely generated torsion-free nilpotent, and residually $p$-finite. Hence, it follows that they are residually $\mathcal{C}$.
\end{enumerate}
\end{proof}
\end{lemma}
In particular, if $\mathcal{C}$ satisfies the Gruenberg condition then $\mathcal{C}$ is closed under finite direct products. This shows that the second condition in the
definition of root class is redundant.

\medskip
Recently, a new characterization of root classes was given.

\begin{theorem}\emph{\cite[Theorem 1]{Sokolov2013}} \label{Sokolov-equivalences}
Let $\mathcal{C}$ be a class of groups closed under taking subgroups, then the following are equivalent:
\begin{enumerate}
 \item $\mathcal{C}$ satisfies the Gruenberg condition (that is, $\mathcal{C}$ is a root class);
 \item $\mathcal{C}$ is closed under unrestricted wreath products;
 \item $\mathcal{C}$ is closed under extensions and for any two groups $G, H \in \mathcal{C}$ the group $\prod_{h\in H}G$ is again in $\mathcal{C}$.
\end{enumerate}
\end{theorem}
\noindent
In the following lemma we use the equivalence between the first and the third condition of this characterization.

\begin{lemma}\label{residuallyroot-by-C}
Let $\mathcal{C}$ be a class of groups that contains a root class $\mathcal{R}$.
Suppose that
\begin{enumerate}
 \item every $\mathcal{R}$-by-$\,\mathcal{C}$ group belongs to $\mathcal{C}$;
 \item for every group in $\mathcal{C}$ there is a group in $\mathcal{R}$ of the same cardinality.
\end{enumerate}
Then a (residually $\mathcal{R}$)-by-$\,\mathcal{C}$ group is residually $\mathcal{C}$.
\begin{proof}
Let $G$ be a (residually $\mathcal{R}$)-by-$\,\mathcal{C}$ group, with a residually $\mathcal{R}$ 
normal subgroup $N$ such that~$G/N\in\mathcal{C}$. 

Consider a non-trivial element $g\in G$. If $g\notin N$ then $g$ is not mapped to the identity of $G/N$ by the canonical projection. 
If $g\in N\setminus\{e\}$, then there 
exists a normal subgroup $K\trianglelefteq N$ such that $g\notin K$ and $N/K\in \mathcal{R}$.
Let $S$ be a set of representatives of $N$ in $G$, then the group $G/N\in\mathcal{C}$ has the same cardinality of $S$. The hypotheses imply that
there exists a group $\Gamma\in \mathcal{R}$ of the same cardinality of $G/H$, so
$$\lvert\Gamma\rvert=\lvert G/N\rvert=\lvert S\rvert.$$
Consider the normal subgroup $L:=\bigcap_{s\in S}K^s$ of $G$, then 
$$N/L\hookrightarrow\prod_{s\in S} N/K^s\cong \prod_{s\in S} N/K.$$
The third condition of Theorem \ref{Sokolov-equivalences} implies that 
$\prod_{s\in S} N/K\in \mathcal{R}$, because $\lvert S\rvert=\lvert \Gamma\rvert$ and $\Gamma\in\mathcal{R}$. As $\mathcal{R}$ 
is closed under taking subgroups, this implies that $N/L\in \mathcal{R}$.

By the Third Isomorphism Theorem for groups $$\frac{G/L}{N/L}\cong G/N,$$ 
so the group $G/L$ is the extension of the group $N/L$, which is in $\mathcal{R}$, by the group $G/N$, which is in $\mathcal{C}$. Thus
$G/L\in \mathcal{C}$.

As $L\leq N$ we have that $g\notin L$. Thus, the surjective homomorphism 
$$\varphi\colon G\twoheadrightarrow G/L,\qquad \varphi(g):=gL \quad \forall g\in G$$
maps $g$ to a non-trivial element of $G/L$, and $G$ is residually $\mathcal{C}$.
\end{proof}
\end{lemma}

\begin{corollary1}
If a group is (residually solvable)-by-amenable then it is residually amenable. 

If a group is (residually solvable)-by-(elementary amenable) then it is residually elementary amenable. 
\begin{proof}
To apply Lemma \ref{residuallyroot-by-C}, we note that both the classes of amenable groups and elementary amenable groups are closed under extensions, 
and they contain the root class of solvable groups.
Moreover, for every cardinal number there exist solvable groups of that cardinality, so the hypotheses 
are satisfied and we can conclude. 
\end{proof}
\end{corollary1}

In particular, it follows that
\begin{corollary}\label{free-by-A}
If a group is free-by-amenable then it is residually amenable.
\begin{proof}
Magnus proved that free groups are residually torsion-free nilpotent \cite{Mag}, hence they are residually solvable. Apply now Corollary \ref{RS-by-A}.
\end{proof}
\end{corollary}

\begin{theorem2}
Let $\mathcal{C}$ be a class of groups that contains a root class $\mathcal{R}$. Suppose that
\begin{enumerate}
 \item $\mathcal{C}$ is closed under finite direct products;
 \item every $\mathcal{R}$-by-$\,\mathcal{C}$ group belongs to $\mathcal{C}$;
 \item for every group in $\mathcal{C}$ there is a group in $\mathcal{R}$ of the same cardinality.
\end{enumerate} 
Then the class of residually $\mathcal{C}$ groups is closed under free products.
\begin{proof}
Let $\{G_i\}_{i\in I}$ be a family of residually $\mathcal{C}$ groups. Without loss of generality we can suppose that $I$ is finite. In fact,
if $I$ is not finite, let $w\in \ast_{i\in I}G_i$ be a non-trivial reduced word and let $J\subseteq I$ be
the finite set of those indices of $I$ that appear in $w$. Consider then the surjective homomorphism
$$\psi\colon \ast_{i\in I}G_i\twoheadrightarrow \ast_{j\in J}G_j$$ 
such that $\psi(g_j)=g_j$ for $j\in J$ and $\psi(g_i)=e$ for $i\in I\setminus J$.
Via $\psi$ we see that if the theorem holds when the index set is finite, it will hold in full generality. 
We prove the theorem for $\lvert I\rvert=2$, as the same argument applies to the other cases. 

Moreover, without loss of generality we can suppose that $G_1$ and $G_2$ are groups in $\mathcal{C}$. In fact,
consider a reduced non-trivial word $w\in G_1\ast G_2$ containing letters from both groups. Let 
$$A_i:=\{g\in G_i\mid g\text{ appears in }w\}, \qquad \text{ for } i=1,2.$$
The sets $A_i$ are finite, non-empty subsets of $G_i$, and they do not contain the identity of $G_i$.
The groups $G_i$ are residually $\mathcal{C}$, and $\mathcal{C}$ is closed under finite direct products. This means that
they are fully residually $\mathcal{C}$, that is, there exist groups $C_1,C_2\in\mathcal{C}$ and surjective homomorphisms 
$$\varphi_1\colon G_1\twoheadrightarrow C_1\qquad \varphi_2\colon G_2\twoheadrightarrow C_2$$ 
such that 
$$\varphi_1(A_1)\subseteq C_1\setminus\{e_{G_1}\},\qquad \varphi_2(A_2)\subseteq C_2\setminus\{e_{G_2}\}$$
and that $\varphi_i\restriction_{A_i}$ is injective.
Define the surjective homomorphism $\varphi\colon G_1\ast G_2\twoheadrightarrow C_1\ast C_2$
such that
$$\varphi\restriction_{G_i}=\varphi_i\qquad\text{ and }\qquad \varphi(g_1\dots g_r)=\varphi(g_1)\dots\varphi(g_r)
\quad\forall g_i\in G_{\iota_i}, \forall r\in\N.$$
This surjective homomorphism $\varphi$ maps $w$ to a non-trivial element of $C_1\ast C_2$, which is a free product 
of groups belonging to $\mathcal{C}$.

Hence suppose without loss of generality that $G_1,G_2\in\mathcal{C}$ and consider the canonical projection 
$$\pi\colon G_1\ast G_2\twoheadrightarrow G_1\times  G_2.$$ 
The class $\mathcal{C}$ is 
closed under finite direct products, so the image of $\pi$ belongs to $\mathcal{C}$.
Let $K=\ker\pi$, then~$K$ is free by the 
Kuro\v{s} Subgroup Theorem \cite[Theorem 6.3.1]{Rob95}, because it has trivial intersection with every $G_i$. 
Free groups are residually $\mathcal{R}$, because
$\mathcal{R}$ is a root class, so it follows that the free product $G_1\ast G_2$ is (residually $\mathcal{R}$)-by-$\,\mathcal{C}$.
Apply Lemma \ref{residuallyroot-by-C} to conclude that $G_1\ast G_2$ is residually $\mathcal{C}$.
 
\end{proof}
\end{theorem2}

\begin{corollary3}
The classes of residually amenable groups and of residually elementary amenable groups are closed under free products.
\begin{proof}
Let $\mathcal{C}$ be the class of amenable groups, or the one of elementary amenable groups, 
$\mathcal{R}$ be the class of solvable groups and apply Theorem \ref{theoremResi}.
\end{proof}
\end{corollary3}

\smallskip

Theorem \ref{theoremResi} guarantees sufficient conditions for a free product of residually $\mathcal{C}$ groups to be again residually $\mathcal{C}$.
In general, this is not always the case. A pathological counterexampleis given by the class of abelian groups.

For a non-trivial example when the conclusion of Theorem \ref{theoremResi} fails, 
\cite[Theorem 6]{Ba67} guarantees that a free product of two non-trivial groups is residually free if and only if the two groups are fully residually
free.
So consider a residually free group that is not fully residually free, for instance $G=F_2\times\Z$. Then
$G\ast G$ is a free product of two residually free groups, and it is not residually free. 
Note that, in both cases, the two classes considered are not closed under extensions.

\medskip\medskip

In \cite{HirThu75} (see also \cite[III. 13.]{dlHarpe}) the following class of groups is considered. Let $\mathcal{C}$ be the smallest class of 
groups satisfying the following conditions:
\begin{itemize}
 \item[$(c1)$] $\mathcal{C}$ contains all amenable groups; 
 \item[$(c2)$] if $G,H\in \mathcal{C}$ then $G\ast H\in\mathcal{C}$;
 \item[$(c3)$] if $G$ is virtually $\mathcal{C}$ then $G\in \mathcal{C}$.
\end{itemize}
It is known that the fundamental group
\begin{equation}\label{surface.groups}
G_k:=\bigl\langle a_1,b_1,\dots, a_k,b_k\mid\prod_{i=1}^r[a_i,b_i]\bigr\rangle
\end{equation}
of a closed orientable surface of genus $k\geq 2$ does not belong to $\mathcal{C}$ \cite{HirThu75}.

Corollary \ref{corollary.resam} and Lemma \ref{virtually.residually.amenable} imply that the class of residually amenable groups 
satisfy the previous conditions $(c2)$ and $(c3)$, and condition $(c1)$ is trivially satisfied. Thus 
$$\mathcal{C}\subsetneq\{\text{residually amenable groups}\},$$
where the inclusion is strict because the groups given by the presentation 
\eqref{surface.groups} are residually finite \cite[\S 1.3]{Ba62}, hence residually amenable.

Moreover
\begin{equation*}\label{Cnotreselam}
\mathcal{C}\nsubseteq\{\text{residually elementary amenable groups}\},
\end{equation*}
in fact this inclusion is not true because there exist finitely generated simple amenable groups \cite{JuMo2013} that are not elementary amenable 
(hence, they are not residually elementary amenable).
These groups are not elementary amenable because they are infinite and \cite[Corollary 2.4]{Chou80} guarantees that a finitely generated simple 
elementary amenable group is finite.

Also the other inclusion does not hold, that is
\begin{equation*}\label{reselamnotC}
\{\text{residually elementary amenable groups}\}\nsubseteq \mathcal{C},
\end{equation*}
because the groups given by \eqref{surface.groups} are residually finite, hence residually elementary amenable.

\medskip

We now prove Theorem \ref{theoremLocalEmbeddings}. Although the class of residually free groups is not closed under taking free products, 
Theorem \ref{theoremLocalEmbeddings} guarantees that a free product of groups that are locally embeddable into free groups is 
again locally embeddable into free groups.

\begin{theorem4}
Let $\mathcal{C}$ be a class of groups, and suppose that
\begin{itemize}
 \item[$(a)$] $\mathcal{C}$ is closed under free products
\end{itemize} 

or
\begin{itemize}
 \item[$(b)$] $\mathcal{C}$ is closed under finite direct products and the free product of residually $\mathcal{C}$ groups is residually $\mathcal{C}$.
\end{itemize}
Then the class of LE-$\mathcal{C}$ groups is closed under free products.
\begin{proof}
Let $\{G_i\}$ be a family of groups locally embeddable into $\mathcal{C}$. As in the proof of Theorem \ref{theoremResi}, we can suppose that
the index set $I$ is finite, and the argument for $\lvert I\rvert=2$ applies to all other cases.
Let $K=\{w_1,\dots, w_r\}$ be a finite set of reduced words of $G_1\ast G_2$, and consider the finite subsets 
$$K_i:=\{e_{G_i}\}\cup\{g\in G_i\mid g \text{ appears in some }w_i \text{ for }j=1,\dots,r\}\subseteq G_i,\qquad \text{ for } i=1,2.$$
The groups $G_i$ are locally embeddable into $\mathcal{C}$, hence
there exist $C_i\in \mathcal{C}$ and $K_i$-almost-homomorphisms $\varphi_i\colon G_i\to C_i$.
Moreover $e_{G_i}\in K_i$, thus $\varphi_i(e_{G_i})=e_{C_i}$. This implies that 
\begin{equation}\label{welldefined}
 \varphi_i(g)\neq e_{C_i}\quad \forall g\in K_i\setminus\{e_{G_i}\}.
\end{equation}
Then the map $\varphi\colon G_1\ast G_2\to C_1\ast C_2$ defined by 
$$\varphi\restriction_{G_i}=\varphi_i\qquad\text{ and }\qquad \varphi(g_1\dots g_s)=\varphi(g_1)\dots\varphi(g_s)
\quad\forall g_i\in G_{\iota_i},\forall s\in\N$$
is well defined in view of Equation \eqref{welldefined}.
Let us see that $\varphi$ is a $K$-almost-homomorphism. 

Consider the reduced words
$$w=g_1\dots g_s\in K,\qquad g_i\in G_{\iota_i}\qquad \text{ and }\qquad w'=h_1\dots h_t\in K,\qquad h_i\in G_{\kappa_i}$$ 
and suppose that $\varphi(w)=\varphi(w')$. 
It follows that $s=t$. Moreover, as $w$ and $w'$ are reduced words, also $\varphi(w)$ and $\varphi(w')$
are reduced. This implies that $\iota_j=\kappa_j$ for all $j=1,\dots,s$ and that
$$\varphi_{\iota_1}(g_1)=\varphi_{\iota_1}(h_1),\quad \dots, \quad\varphi_{\iota_s}(g_s)=\varphi_{\iota_s}(h_s).$$
As the maps $\varphi_{\iota_j}$ are $K_j$-almost-homomorphisms and $g_j,h_j\in K_{\iota_j}$, it follows that $g_j=h_j$ for all $j=1,\dots, s=t$.
This means that $w=w'$, that is, $\varphi\restriction_K$ is injective.

Moreover, condition $(a1)$ of the definition of $K$-almost-homomorphism is satisfied, because the maps $\varphi_j$ satisfy 
the same conditions for $K_j$, respectively. 
Thus $\varphi$ is a $K$-almost-homomorphism.

\medskip

If condition $(a)$ holds, then $C_1\ast C_2\in\mathcal{C}$ and so $G_1\ast G_2$ is locally embeddable into $\mathcal{C}$.
On the other hand, if condition $(b)$ holds then $C_1\ast C_2$ is residually $\mathcal{C}$. 
As the class $\mathcal{C}$ is closed under finite direct products, the group $C_1\ast C_2$ is fully residually $\mathcal{C}$, that is,
there exists a surjective homomorphism 
$\psi\colon C_1\ast C_2\twoheadrightarrow D\in\mathcal{C}$
which is injective on the finite subset $\varphi(K)$ of $C_1\ast C_2$.
Thus, the composition $\psi\circ\varphi\colon G_1\ast G_2\to D$ is a $K$-almost-homomorphism, and $G_1\ast G_2$ 
is LE-$\mathcal{C}$.
 
\end{proof}
\end{theorem4}

Conditions $(a)$ and $(b)$ of the previous theorem are independent one from the other. In fact the class of free groups satisfies the first one 
but not the second, and the class of finite (or amenable, etc.) groups satisfies the second but not the first.
As for the residual case, a class of groups $\mathcal{C}$ such that being locally embeddable into $\mathcal{C}$ is not preserved under
free products is the class of abelian groups.

Corollary \ref{LEA-free} follows immediately from Theorem \ref{theoremLocalEmbeddings}.

\section{The pro-$\mathcal{C}$ topology and applications to HNN extensions and amalgamated free products}\label{section-proamenable}

In this section, we consider the pro-$\mathcal{C}$ topology and the pro-$\mathcal{C}$ completion of a group $G$, we investigate some properties 
of these constructions for residually $\mathcal{C}$ groups and we give applications to HNN extensions and amalgamated free products. 

In \cite{ArHKS} the analogous concept of true prosoluble completion for a group is studied, and in \cite{RibZal} 
pro-$\mathcal{C}$ groups are considered, where $\mathcal{C}$ are various classes of finite groups. We give reference to these works and to
the references therein for more about pro-$\mathcal{C}$ completions.

\smallskip
Let $\mathcal{C}$ be a class of groups, define the \emph{pro}-$\mathcal{C}$ \emph{topology} of a group $G$ as the topology with the following basis at $e_G$
\begin{equation*}\label{proamenable-def}
\mathcal{B}:=\bigl\{N\trianglelefteq G\mid G/N\in\mathcal{C}\bigr\}. 
\end{equation*}
The pro-$\mathcal{C}$ topology of a group is a linear topology, hence it is a group topology, namely, the group operations
$$\mu\colon G\times G\to G,\qquad \iota\colon G\to G$$
defined by $\mu(g,h):=gh$ and $\iota(g):=g^{-1}$ are continuous maps, where $G\times G$ is endowed with the product topology.

If $\mathcal{C}$ is closed under subgroups and finite direct products, then the family $\mathcal{B}$ is directed with respect to inclusion:
if $N_1$ and $N_2$ are elements of $\mathcal{B}$ then also $N_1\cap N_2$ is.

Define hence the \emph{pro}-$\mathcal{C}$  \emph{completion} of a group $G$ as the inverse limit
\begin{equation*}\label{proamenable.completion}
\widehat{G}:=\varprojlim_{N\in\mathcal{B}} G/N
\end{equation*}
with respect to the canonical surjective homomorphisms $G/N\twoheadrightarrow G/M$, for $N\leq M$ in $\mathcal{B}$,
and let $\eta$ be the natural homomorphism 
$$\eta\colon G\to \widehat{G},\qquad\eta(g):=\bigl(gN\bigr)_{N\in\mathcal{B}}.$$
We have the following characterizations. Let $G$ be a group and $\mathcal{C}$ be a class closed under subgroups and 
finite direct products, then the following conditions are equivalent:
\begin{itemize}
 \item[$(c1)$] $G\in\mathcal{C}$;
 \item[$(c2)$] $\{e_G\}\in\mathcal{B}$, that is, the pro-$\mathcal{C}$ topology over $G$ is discrete;
 \item[$(c3)$] $\widehat{G}=G$ and $\eta\colon G\to G$ is the identity map.
\end{itemize}
Also the following conditions are equivalent:
\begin{itemize}
 \item[$(rc1)$] $G$ is residually $\mathcal{C}$;
 \item[$(rc2)$] the pro-$\mathcal{C}$ topology over $G$ is Hausdorff;
 \item[$(rc3)$] $\eta\colon G\to \widehat{G}$ is injective.
\end{itemize}

For the rest of the section, also if not explicitly stated, $\mathcal{C}$ is assumed being closed under subgroups and finite direct products.
The proof of the following lemma carries over from the proof of \cite[Theorem 3.1]{Hall}.

\begin{lemma}\label{Hallreference}
Let $G$ be a residually $\mathcal{C}$ group and consider its 
pro-$\mathcal{C}$ topology. A subgroup $H$ contains an open set if and only if it contains some subgroup $N\in\mathcal{B}$.

Moreover, in such a case, $H$ is both open and closed.
\end{lemma}

We underline the following immediate corollary:

\begin{corollary}\label{cor.helpful}
Let $G$ be a residually $\mathcal{C}$ group, $H$ a subgroup and $N\in\mathcal{B}$. Then the subgroup $HN$ is both open and closed in the pro-$\mathcal{C}$ topology.
\end{corollary}

For a subgroup $H\leq G$ the following three equivalent properties characterize to be closed in the pro-$\mathcal{C}$ topology.
\begin{enumerate}
 \item for all $g\in G\setminus H$ there exists a surjective homomorphism $\varphi\colon G\twoheadrightarrow C$ 
 such that $C\in\mathcal{C}$ and $\varphi(g)\notin \varphi(H)$;
 \item for all $g\in G\setminus H$ there exists $N\in\mathcal{B}$ such that $g\notin HN$;
 \item $H=\bigcap \bigl\{HN\mid N\in\mathcal{B}\bigr\}$.
\end{enumerate}

When $\mathcal{C}=\{\text{finite groups}\}$, a subgroup satisfying one of the above conditions is sometimes said to be separable in the profinite topology.

\medskip\medskip

We now turn our attention to HNN extensions and amalgamated free products. 
Recall that, given a group $A$, two subgroups $H$ and $K$ and $\varphi\colon H\to K$ an isomorphism,
the HNN extension with base group $A$ and amalgamated subgroups $H$ and $K$ is the group given by the presentation
$$A\ast_\varphi:=\langle A,t\mid t^{-1}ht=\varphi(h)\quad \forall h\in H\rangle.$$
We call the HNN extension \emph{special} if $H=K$ and $\varphi=\id_H$.

Theorem \ref{specialHNN} extends
\cite[Theorem 4.2]{Tieudjo2010} via Lemma \ref{residuallyroot-by-C} to classes of groups $\mathcal{C}$ that satisfy the same hypotheses
of that lemma.

\begin{theorem5}
Let $A$ be a group, $H\leq A$ and consider the special HNN extension $G:=A\ast_{\id}$ amalgamating the subgroup $H$. 
Let $\mathcal{C}$ be a class of groups that contains a root class $\mathcal{R}$. Suppose that
\begin{enumerate}
 \item $\mathcal{C}$ is closed under subgroups and finite direct products;
 \item every $\mathcal{R}$-by-$\,\mathcal{C}$ group belongs to $\mathcal{C}$;
 \item for every group in $\mathcal{C}$ there is a group in $\mathcal{R}$ of the same cardinality.
\end{enumerate}
Then $G$ is residually $\mathcal{C}$ if and only if $A$ is residually $\mathcal{C}$ and $H$ is closed in the pro-$\mathcal{C}$ topology of $A$.
\begin{proof}
We first assume that $A\in\mathcal{C}$ and we prove that the HNN extension $G$ is residually $\mathcal{C}$. 
In this case the pro-$\mathcal{C}$ topology of $A$ is clearly discrete, hence
all the subgroups are closed.

Consider the surjective homomorphism $\pi\colon G\twoheadrightarrow A$ defined by $\pi(a)=a$ for all $a\in A$ and $\pi(t)=e_A$. 
This is a homomorphism because the considered HNN extension is special.
Let $K:=\ker \pi$, then $A\cap K=\{e_A\}$, and for all $g\in G$ we have that
$$g^{-1}Ag\cap K=g^{-1}Ag\cap g^{-1}Kg=g^{-1}(A\cap K)g=\{e_A\}.$$
Hence, $K$ is a free group because it acts freely on the Bass-Serre tree of the HNN extension $A\ast_\id$. 
This means that $G$ is free-by-$\,\mathcal{C}$ and in particular (residually~$\mathcal{R}$)-by-$\,\mathcal{C}$. Lemma \ref{residuallyroot-by-C} 
implies that the group $G$ is residually $\mathcal{C}$.

\medskip\medskip

Suppose that the base group $A$ is residually $\mathcal{C}$ and that $H$ is closed in $A$, we want to show that $G$ is residually~$\mathcal{C}$.
Let $g$ be a non-trivial element of $G$. By Britton's lemma $g$ has the following reduced form
\begin{equation}\label{reducedform}
g=a_0t^{\varepsilon_1}a_1\dots t^{\varepsilon_n}a_n,
\end{equation}
where
\begin{itemize}
 \item $n$ is minimal, $a_i\in A$ and $\varepsilon_i=\pm1$ for all $i=0,\dots,n$;
 \item $a_i\in \bigl(A\setminus H\bigr)\cup \{e\}$ for all $i=1,\dots,n$;
 \item if $\varepsilon_{i}\varepsilon_{i+1}=-1$ for some $i=1,\dots,n-1$ then $a_i\neq e_A$ (and in particular $a_i\notin H$).
\end{itemize}

If $g=t^\varepsilon$ for some $\varepsilon \in\Z\setminus\{0\}$, then for all normal subgroups $N\trianglelefteq A$ the element $g$ is mapped 
non-trivially by $\pi\colon G\twoheadrightarrow G_N$ to the special HNN extension 
$$G_N:=A/N\ast_{\bar{\id}}=\langle A/N,t\mid t^{-1}xt=x,\,\forall x\in HN/N\rangle,$$
where $\bar{\id}\colon HN/N\to HN/N$ is the identity map.
Let $N\in\mathcal{B}$, then $G_N$ is residually $\mathcal{C}$ because the base group $A/N$ is in $\mathcal{C}$.
Hence there exists a surjective homomorphism $\varphi\colon G_N\twoheadrightarrow B$, where $B\in\mathcal{C}$ such that $\varphi(t^\varepsilon)\neq e_B$. 
Thus, the composition of $\varphi\circ\pi$ maps $g$ to a non-trivial element of the group $B\in\mathcal{C}$.

If $g\neq t^\varepsilon$, then in its reduced form given by \eqref{reducedform} there are some $a_i\neq e_A$. As $H$ is closed in $A$ and hence it is separable in the pro-$\mathcal{C}$ topology,
there exists a surjective homomorphism $\varphi\colon A\twoheadrightarrow B$, with $B\in\mathcal{C}$, such that for every non-trivial $a_i$ (where the elements 
$a_i$ are the elements in the reduced
form of $g$ of Equation~\eqref{reducedform}) 
\begin{equation}\label{nontrivialcondition}
\varphi(a_i)\notin \varphi(H)=HN/N,
\end{equation}
where $N:=\ker\varphi$. Extend $\varphi\colon A\twoheadrightarrow B$ to $\bar{\varphi}\colon A\ast_{\id}\twoheadrightarrow B\ast_{\id}$, where in $B\ast_{\id}$ the identity map is
$\id\colon \varphi(H)\to\varphi(H)$, that is, the identity map of $HN/N$. The element $\bar{\varphi}(g)$ is non-trivial in $B\ast_{\id}$, because of
the conditions in Equation~\eqref{nontrivialcondition}, and $B\ast_{\id}$ is residually $\mathcal{C}$ by the first part of the theorem, because 
$B\in\mathcal{C}$.

This means that $G$ is residually (residually $\mathcal{C}$), that is, $G$ is residually $\mathcal{C}$.
\medskip\medskip

To conclude, suppose that the special HNN extension $G=A\ast_{\id}$ is residually $\mathcal{C}$, we want to show that $A$ is residually
$\mathcal{C}$ and that the subgroup $H$ is closed in the pro-$\mathcal{C}$ topology of $A$.

As $A\leq G$, it is residually $\mathcal{C}$.
Suppose that $H$ is not closed in $A$, so there exists an element $\tilde{a}\in A\setminus H$ such that $\tilde{a}\in HN$ for all the normal subgroups 
$N\in\mathcal{B}$.
Consider the non-trivial element $g=t^{-1}\tilde{a}t\tilde{a}^{-1}\in G$.
The group $G$ is residually $\mathcal{C}$, so there exists a normal subgroup $K\trianglelefteq G$ such that $g\notin K$ and $G/K\in\mathcal{C}$.

Let $N:=K\cap A\trianglelefteq A$, then $A/N\in\mathcal{C}$, because
$$A/N=\frac{A}{A\cap K}\cong \frac{AK}{K}\leq G/K.$$
Extend the canonical projection $\pi\colon A\twoheadrightarrow A/N$ to 
$$\bar{\pi}\colon A\ast_{\id}\twoheadrightarrow \bigl(A/N\bigr)\ast_{\id},$$
where in the second HNN extension the identity map is $\id\colon HN/N\to HN/N$. 
As $A/N\in\mathcal{C}$, it follows that $\tilde{a}\in HN$. So there exist $h\in H$
and $n\in N$ such that $\tilde{a}=hn$.
Thus
$$\bar{\pi}(\tilde{a})=\pi(\tilde{a})=\pi(h)\pi(n)=\pi(h)\in\pi(H)=HN/N$$
and hence
$$\bar{\pi}(g)=\bar{\pi}(t^{-1}\tilde{a}t\tilde{a}^{-1})=\underbrace{t^{-1}\pi(h)t}_{=\pi(h)}\pi(h^{-1})=\pi(h)\pi(h^{-1})=e.$$
Thus $g\in \ker\pi_N\leq K$ and in particular $g\in K$. 
This is a contradiction with the initial choice of the normal subgroup $K$, and this contradiction arises from supposing that $H$ is not
closed in $A$.

\end{proof}
\end{theorem5}

\begin{corollary}
Let $\mathcal{C}$ be a class of groups as in the previous theorem, $A$ a residually $\mathcal{C}$ group and $H\leq A$ 
a finite subgroup of $A$. Then the special HNN extension $A\ast_{\id}$ with associated subgroup $H$ is residually $\mathcal{C}$.

\end{corollary}

The proof of Corollary \ref{cor-HNNspecial} follows from Theorem \ref{specialHNN}. 

\medskip\medskip

Not surprisingly, a similar statement is true for amalgamated free products. Let $H\leq A$, $K\leq B$ be groups and $\varphi\colon H\to K$ an isomorphism.
We call the amalgamated free product $A\ast_{\varphi}B$ a \emph{double} if $A=B$ and $\varphi=\id_H$. In this case, we write $A\ast_HA$ for $A\ast_\varphi A$, and we denote by $\overline{A}$
the right hand side copy of $A$ in the double.

\begin{theorem6}
Let $A$ be a group, $H\leq A$ a subgroup and consider the double $G:=A\ast_H A$. 
Let $\mathcal{C}$ be a class of groups that contains a root class $\mathcal{R}$. Suppose that
\begin{enumerate}
 \item $\mathcal{C}$ is closed under subgroups and finite direct products;
 \item every $\mathcal{R}$-by-$\,\mathcal{C}$ group belongs to $\mathcal{C}$;
 \item for every group in $\mathcal{C}$ there is a group in $\mathcal{R}$ of the same cardinality.
\end{enumerate}
Then $G$ is residually $\mathcal{C}$ if and only if $A$ is residually $\mathcal{C}$ and $H$ is closed in the pro-$\mathcal{C}$ topology of $A$.
\begin{proof}
Suppose first that $A\in\mathcal{C}$. Let $\bar{\id}\colon\overline{A}\to A$ be the isomorphism sending $\bar{a}$ to $a$.

By the universal property of amalgamated free products, there exists a unique homomorphism $\varphi\colon G\to A$ such that
$\varphi\restriction_A=\id$, $\varphi\restriction_{\overline{A}}=\bar{\id}$:
$$\xymatrix{
A\ar[dr]_{\id}\ar@{^{(}->}[r]&G\ar@{.>}[d]^{\varphi}&\overline{A}\ar[dl]^{\bar{\id}}\ar@{_{(}->}[l]\\
&A&}$$
Let $K=\ker\varphi$, then $K$ is generated by the elements $a\bar{a}^{-1}$, for all $a\in A$. 
Moreover $K\cap A=K\cap \overline{A}=\{e\}$ and hence $K$ is free. 
The homomorphism $\varphi$ is surjective, so it follows that $G$ is free-by-amenable and thus residually amenable by Corollary \ref{free-by-A}.

\medskip\medskip

Suppose that $A$ is a residually $\mathcal{C}$ group and that $H\leq A$ is closed in the pro-$\mathcal{C}$ topology of $A$.
Let $g\in G$ be a non-trivial element of the double. If $g\in A$ then $\varphi(g)=g\in A\setminus\{e_A\}$. 
The group $A$ is residually $\mathcal{C}$, so there exists a group $B\in\mathcal{C}$ and a surjective homomorphism $\psi\colon A\to B$ 
such that $\psi(g)\neq e_B$. The same argument works for $g\in\overline{A}$.

If $g\notin A\cup \overline{A}$ then
$$g=a_1b_1\dots a_nb_n \qquad a_i\in A\setminus H,\, b_i\in \overline{A}\setminus\overline{H}.$$
The subgroup $H$ is closed in the pro-$\mathcal{C}$ topology of $A$, thus there exist a group $C\in\mathcal{C}$ and a surjective homomorphism 
$\theta\colon A\to C$ such that 
$$\theta(a_i)\notin\theta(H),\qquad\bar{\theta}(\bar{a_i})\notin\bar{\theta}(\overline{H}),\qquad\forall i=1,\dots,n.$$
Consider the surjective homomorphism
$$\Theta\colon G\twoheadrightarrow C\ast_{\theta(H)}C$$
such that
$$g\mapsto\Theta(g)=\theta(a_1)\bar{\theta}(b_1)\dots\theta(a_n)\bar{\theta}(b_n).$$
As $\theta(a_i)\notin \theta(H)$ and $\bar{\theta}(b_i)\notin\bar{\theta}(\overline{H})$, it follows that $\Theta(g)$ is not trivial in
$C\ast_{\theta(H)}C$.
By the first part of this proof $C\ast_{\theta(H)}C$ is residually $\mathcal{C}$, hence there exists a quotient $B\in\mathcal{C}$ of $G$ where
$g$ is mapped non-trivially.

\medskip\medskip

Suppose now that the amalgamated free product $G$ is residually $\mathcal{C}$ and thus that $A$ is residually $\mathcal{C}$.
We need to prove that the amalgamated subgroup $H$ is closed in the pro-$\mathcal{C}$ topology. Suppose it is not, then there exists an element 
$a\in A\setminus H$ such that $a\in HN$ for all normal subgroups $N\trianglelefteq A$ with $A/N\in\mathcal{C}$.
Consider the element $[a,\bar{a}]\in G\setminus (A\cup \overline{A})$.

As the group $G$ is residually $\mathcal{C}$, there exists a normal subgroup $K\trianglelefteq G$ such that $[a,\bar{a}]\notin K$ and 
$G/K\in\mathcal{C}$.
Let $N:=K\cap A$, then $N\trianglelefteq A$ and the quotient 
$$A/N\cong \frac{AK}{K}\leq \frac{G}{K}\in\mathcal{C}.$$
Thus $a\in HN$, so there exist $h\in H$ and $k\in N$ such that $a=hk$. Moreover $N\trianglelefteq A$, so there exist elements
$n,m\in N$ such that $h^{-1}\bar{k}^{-1}=nh^{-1}$ and $k\bar{h}=\bar{h}m$. This means that
\begin{equation*}
\begin{split}
[a,\bar{a}]&=(hk)^{-1}(\bar{h}\bar{k})^{-1}hk\bar{h}\bar{k}=k^{-1}h^{-1}\bar{k}^{-1}\bar{h}^{-1}hk\bar{h}\bar{k}\\
&=k^{-1}(h^{-1}\bar{k}^{-1})(k\bar{h})\bar{k}=
k^{-1}(nh^{-1})(\bar{h}m)\bar{k}=k^{-1}nm\bar{k}\in N=K\cap A\leq K.
\end{split} 
\end{equation*}
This is a contradiction, because $[a,\bar{a}]$ survives in the quotient $G/K$. Hence the subgroup $H$ is closed in the
pro-$\mathcal{C}$ topology of $A$, and the proof is completed.
\end{proof}
\end{theorem6}

\begin{corollary}\label{this}
Let $\mathcal{C}$ be a class of groups as in the previous theorem, $A$ a residually $\mathcal{C}$ group and $H\leq A$ 
a finite subgroup of $A$. Then the double $A\ast_HA$ is residually $\mathcal{C}$.

\end{corollary}

Corollary \ref{cor.double} follows from Theorem \ref{theo.bingobongo}.

\begin{remark}
Theorem \ref{theo.bingobongo} is no longer valid if the amalgamated free product is not a double. In fact, in \cite{BuMo} infinite simple groups
are constructed as amalgamated free products of free groups (hence, residually finite groups) amalgamating subgroups of finite index (hence,
closed in the profinite / proamenable topology of the given free groups). 

Thus an amalgamated free product of residually amenable groups, amalgamating two closed subgroups, may be not residually amenable in general.
\end{remark}

\medskip\medskip

Consider a solvable group $A$ with two isomorphic subgroups $\varphi\colon H\to K$, and suppose that $\varphi$ is the restriction of an 
automorphism of $A$.
It was proved that, in these hypotheses, the HNN extension $A\ast_\varphi$ is residually solvable \cite[Lemma 1.2]{RaVa}.
Proposition \ref{stabilityHNN} is the generalization of this fact just mentioned.

\begin{lemma6}
Let $A$ be a group, $H$, $K$ two subgroups of $A$ and $\varphi\colon H\to K$ an isomorphism. Suppose that there exists an automorphism 
$\alpha\colon A\to A$ such that $\alpha\restriction_H=\varphi$.
If $A$ is amenable then the HNN extension $A\ast_{\varphi}$ is residually amenable. 

\begin{proof}

Consider the two HNN extensions
$$G=\langle A,t\mid t^{-1}ht=\varphi(h)\quad\forall h\in H\rangle,\qquad G^{\star}=\langle A,t\mid t^{-1}at=\alpha(a)\quad\forall a\in A\rangle.$$
As $H\subseteq A$ and $\alpha\restriction_H=\varphi$, the map $\rho\colon G\to G^{\star}$ defined by $\rho(a)=a$ for all $a\in A$ and $\rho(t)=t$ is a well-defined homomorphism. 
Moreover, $G^\star$ is isomorphic to the semidirect product $A\rtimes_\alpha\langle t\rangle$ and thus it is amenable-by-solvable, so amenable.

Let $K=\ker\rho$, we have that $A\cap K=\{e_A\}$, thus $K$ intersects trivially each conjugate of $A$ in $G$ because $K$ is normal in $G$. 
Hence $K$ is a free group, because it acts freely on the Bass-Serre tree of the HNN extension $A\ast_\varphi$.
This means that $G$ is free-by-amenable, so it is residually amenable by Corollary~\ref{free-by-A}.
\end{proof}
\end{lemma6}

With the same hypotheses of Lemma \ref{stabilityHNN}, if the base group is residually amenable then the HNN extension need not 
be residually amenable: Example \ref{ex-amenableNOGRUENBERG} gives the counterexample.
\section*{Acknowledgements}

The author is very grateful to his advisor, Goulnara Arzhantseva, for suggesting the topic, for many motivating questions, valuable remarks
and corrections on previous versions of this paper. 
He wants to thank his colleagues of the Geometric Group Theory Research Group at the University of Vienna 
and Simone Virili for helpful comments and many discussions. Moreover, he thanks Aditi Kar and Nikolay Nikolov for pointing out an error in a previous version of the paper.

The author is supported by the European Research Council (ERC) grant of Prof. Goulnara Arzhantseva, grant agreement no.~259527.

\end{document}